\documentclass[11pt]{article}
\usepackage{amsmath}
\usepackage{amssymb}
\usepackage{amsthm}
\usepackage{verbatim}
\usepackage{enumerate}
\usepackage{listings}
\font\elevenss=cmss11

\font\eightss=cmss8

\font\sixss=cmss8 at 6pt

\newfam\ssfam
\textfont\ssfam=\elevenss \scriptfont\ssfam=\eightss
 \scriptscriptfont\ssfam=\sixss

\theoremstyle{plain}
\newtheorem{thm}{Theorem}

\newtheorem{lem}[thm]{Lemma}

\newtheorem{cor}[thm]{Corollary}

\newtheorem{conj}[thm]{Conjecture}

\theoremstyle{remark}

\def\cS{{\mathcal S}}
\def\Var{{\rm Var} \,}
\def\ee{\varepsilon}
\def\E{{\mathbb E}}
\def\P{{\mathbb P}}
\def\Cox{\hfill \Box}

\def\disp{\displaystyle}
\def\one{{\bf 1}}
\def\|{{\, | \, }}

\def\yy{{\bf y}}
\def\zz{{\bf z}}
\def\Fhat{\hat{F}}
\def\diseq{\stackrel{\cal D}=}

\def\.{\hskip.06cm}
\def\ts{\hskip.03cm}

\begin{document}
\begin{center}
{\large \bf On the longest $k$-alternating subsequence}
\\[5ex]
\end{center}

\begin{center}
Igor Pak\footnote{Department of Mathematics, UCLA, Los Angeles, CA 90095,
\ts\texttt{\{pak\}@math.ucla.edu}}$^,$\footnote{Research supported in part by NSF grant \# DMS-1001842}
and Robin Pemantle\footnote{Department of Mathematics,
University of Pennsylvania, 209 South 33rd Street,
Philadelphia, PA 19104,  \\
\texttt{pemantle@math.upenn.edu}}$^,$\footnote{Research
supported in part by NSF grant \# DMS-1209117}
\end{center}


\noindent{\bf Abstract:}
We show that the longest $k$-alternating substring of a random
permutation has length asymptotic to $2 (n-k) / 3$.





\vskip1.2cm

\subsection*{Introduction}
%
An {\em alternating permutation} is a permutation $\pi \in \cS_n$
satisfying $\pi (1) < \pi (2)  > \pi (3) < \pi (4) > \cdots$.
Alternating permutations have been well studied and enumerated
(see e.g.~\cite{Stan2010}).
Let $L_n$ be the length of the \emph{longest alternating subsequence} 
of a permutation chosen at random uniformly from $\cS_n$.  Motivated by 
the study of longest increasing subsequences, Stanley computed the 
mean and variance of~$L_n$~:
\begin{eqnarray}
\E L_n & = & \frac{4n+1}{6} \\[1ex]  \label{eq:mean}
\Var L_n & = & \frac{8n}{45} - \frac{13}{180} \label{eq:variance}
\end{eqnarray}
for all $n\ge 4$~\cite{Stan2008} (see also~\cite{Stan2007,Romik}).  
In fact, the distribution is asymptotically normal with these 
parameters~\cite{Widom} (see also~\cite[Theorem~2.1]{HoRe2010}).

A {\em $k$-alternating permutation} is a permutation $\pi \in \cS_n$
such that $(-1)^j (\pi (j) - \pi (j+1)) \geq k$ for all $j \in \{1,\ldots,n-1\}$.
In other words, $\pi$ must be alternating and its jumps
$|\pi(j+1) - \pi (j)|$ must all be at least~$k$.
For $k=1$ we get the ordinary alternating permutations.
We learned of $k$-alternating permutations from
D.~Armstrong~\cite{ober2014}, who attributes the
definition to R.~Chen (personal communication, inspired by
a 2011 talk by R.~Stanley).

Let $\pi$ be a uniformly chosen random permutation in $\cS_n$ and let
$L_{n,k} = L_{n,k} (\pi)$ denote the length of the longest
$k$-alternating subsequence of $\pi$.  Armstrong~\cite{ober2014}
made the following conjecture, and verified it via exact computation
for certain small values of $n$ and $k$.
\begin{conj}[Armstrong, 2014]
For all \ts $n \geq 2$ \ts and \ts $k \in \{1,\ldots, n-1\}$, we have:
\begin{equation}
\E L_{n,k} = \frac{4 (n-k) + 5}{6}  \, .
\end{equation}
\end{conj}
In this note we use probabilistic methods to prove the
following asymptotic version of the conjecture.
\begin{thm} \label{th:main}
$$\E L_{n,k} = \frac{2(n-k)}{3} + O(n^{2/3}) \, .$$
\end{thm}

This is proved via the related notion of $x$-alternation for
$x \in (0,1)$, cf.~\cite{ACSS}.  A vector $\yy  = (y_1 , \ldots , y_n) \in [0,1]^n$
is called $x$-alternating if $(-1)^n (y_j - y_{j+1}) \geq x$ for
all $1 \leq j \leq n-1$.  Let $\mu$ denote product Lebesgue
measure on $[0,1]^n$.  Let $\Psi$ be the map taking $\yy \in [0,1]^n$
to the element $\pi \in \cS_n$ defined by
$$\pi (j) = \# \{ i \leq n \, : \, y_i \leq y_j \} \, .$$
A well known fact attributed to R{\'e}nyi (see~\cite{Resn1999}) says
that if $\yy$ has law $\mu$ then $\Psi (\yy)$
is uniformly distributed on $\cS_n$.
Let $L_{n,x} (\yy)$
denote the length of the longest $x$-alternating subsequence
of~$\yy$.  No confusion can result between this and the definition
of $L_{n,k}$ above, provided that we restrict $x$ to $[0,1)$ and
$k$ to positive integral values.  Theorem~\ref{th:main} is a
consequence of the following results.

\begin{lem} \label{lem:select}
Let $Z$ be a binomial random variable with parameters $n$ and $1-x$.  Then
$$L_{n,x} (\yy) \diseq L_{Z,1} \, .$$
In other words, the law of the longest $x$-alternating subsequence
may be exactly simulated by choosing $Z \sim {\rm Bin}(n,1-x)$,
choosing $\pi$ uniformly on $\cS_Z$, and taking the longest alternating
subsequence of $\pi$.
\end{lem}

\begin{cor} \label{cor:select}
\begin{eqnarray}
\E L_{n,x} & = & \frac{2}{3} \ts n(1-x) + \frac{1}{6} \label{eq:mean2} \\[1ex]
\Var L_{n,x} & = & (1-x)(2+5x) \ts \frac{4n}{45} \label{eq:variance2}
\end{eqnarray}
\end{cor}

\noindent{\sc Proof:} Taking expectations in~\eqref{eq:mean}
gives $\E L_{n,x} = (2/3) \E Z + 1/6$, proving~\eqref{eq:mean2}.
The identity $\Var (Y) = E \Var (Y \| Z) + \Var \E (Y \| Z)$
applied to $Y = L_{n,x}$ gives
\begin{eqnarray*}
\Var (L_{n,x}) & = & \E \ts \frac{8Z}{45} - \frac{13}{180} +
   \Var \left ( \frac{2}{3} Z + \frac{1}{6} \right ) \\[2ex]
& = & \frac{8n(1-x)}{45} - \frac{13}{180} + \frac{4}{9} n x (1-x) \\[1ex]
& = & \frac{8 n (1-x) + 20 n x (1-x)}{45}
\end{eqnarray*}
and proves the corollary.   $\Cox$

\begin{lem} \label{lem:compare}
Let $\yy$ be random with law $\mu$.  Denote
\begin{eqnarray*}
x_1 (k,n) & := & k/n - n^{-1/3} \\
x_2 (k,n) & := & k/n + n^{-1/3}
\end{eqnarray*}
Then the following two implications hold with probability $1 - o(1)$
as $n \to \infty$.
\begin{enumerate}[(i)]
\item For all subsequences $\yy'$ of $\yy$, if $\yy'$ is
$x_2$-alternating then $\pi' := \Psi (\yy')$
is $k$-alternating.
\item For all subsequences $\yy'$ of $\yy$, if $\yy'$ is
not $x_1$-alternating then $\pi' := \Psi (\yy')$
is not $k$-alternating.
\end{enumerate}
Consequently, with probability $1-o(1)$,
\begin{equation} \label{eq:sandwich}
L_{n,x_2} (\yy) \leq L_{n,k} (\Psi (\yy)) \leq
   L_{n,x_1} (\yy) \, .
\end{equation}
\end{lem}

\smallskip

\noindent
{\bf Proof of Theorem~\ref{th:main}.} \ts
The theorem follows from
Corollary~\ref{cor:select} and Lemma~\ref{lem:compare}.
Taking expectations in~\eqref{eq:sandwich} we find that
$$\E L_{n,x_2} \leq \E L_{n,k} \leq \E L_{n,x_1} \, .$$
Corollary~\ref{cor:select} then sandwiches $\E L_{n,k}$
between two quantities both of which are asymptotic
to $(2/3) (n-k)$:
\begin{eqnarray*}
\E L_{n,x_j} & = & \frac{2}{3} n (1-x_j) + \frac{1}{6} \\[2ex]
& = & \frac{2}{3} (n-k) + O(n^{2/3}) \,,
\end{eqnarray*}
where $j\in\{1,2\}$.
$\Cox$

\subsection*{Proof of Lemma~\protect{\ref{lem:compare}}}

Let $\Fhat$ denote the empirical distribution
of the values of $\yy$: $\Fhat (t) := n^{-1} \sum_j \one_{y_j \leq t}$.
If~$(i)$ fails then there are $i,j \leq n$ with $|y_i - y_j| \geq x_2$
and $|\pi (i) - \pi (j)| < k$, where $\pi = \Psi (\yy)$.  Letting
$t$ denote the minimum of $y_i$ and $y_j$, this implies that
$\Fhat (t + x_2) - \Fhat (t) < k/n$.  Because
$$\Fhat (t + x_2) - \Fhat (t) =
   \left(\Fhat(t+x_2) - (t+x_2)\right) - \left(\Fhat (t) - t\right) + x_2$$
it follows that
$$|\Fhat (s) - s| > \frac{1}{2} \left(x_2 - \frac{k}{n}\right)
   = \frac{1}{2} \ts n^{-1/3}$$
either for $s = t$ or $s = t+x_2$.
Similarly, if~$(ii)$ fails then there are $i,j \leq n$
with $|y_i - y_j| < x_1$ and $|\pi (i) - \pi (j)| \geq k$,
leading to
$$|\Fhat (s) - s| >  \frac{1}{2} \left(\frac{k}{n} - x_1\right)
   = \frac{1}{2}\ts n^{-1/3}$$
for some $s \in (0,1)$.  In either case,
$$\sup_{s \in [0,1]} |\Fhat (s) - s| > \frac{1}{2} \ts n^{-1/3} \, .$$
But $\sqrt{n} \sup_{s \in [0,1]} |\Fhat (s) - s|$ converges in
distribution to the Kolmogorov-Smirnov statistic (the law of the
maximum of a Brownian bridge).  Because
$n^{-1/3} / n^{-1/2} \to \infty$, this implies that
$$\P \left ( \sup_{s \in [0,1]} |\Fhat (s) - s| > \frac{1}{2} \ts n^{-1/3}
   \right )  \to 0 \, $$
proving the lemma.
$\Cox$

\subsection*{Proof of Lemma~\protect{\ref{lem:select}}}

We begin with another well known fact, attributed to
M.~B{\'o}na in~\cite{Stan2007}: for $\pi \in \cS_n$,
one alternating subsequence $(\pi (i) : i \in A)$ of maximal length
is obtained by selecting $i \in A$ if and only if $i$ is a peak or
a valley, that is,
$\pi (i-1) < \pi (i) > \pi (i+1)$ or $\pi (i-1) > i < \pi (i+1)$,
except that we select $1$ if and only if $\pi (1) < \pi (2)$ 
(see the proof in~\cite[$\S$2]{HoRe2010}).
This generalizes to $k$-alternating subsequences via the following
algorithm which selects the index set $A$ of a $k$-alternating
subsequence of a given permutation $s \in \cS_n$.

{\small
\begin{verbatim}
GREEDY PROVISIONAL ACCEPTANCE:

     Initialize i := 1, j := 2, state := up, A := empty.
     While j <= n do:
        IF (state = up) and s(i) < s(j) < s(j) + k THEN j := j+1 ELSE
        IF (state = up) and s(i) > s(j) THEN i := j , j := j+1 ELSE
        IF (state = up) THEN
             A := A union {i}, i := j, j := j+1, state := down ELSE
        IF s(i) > s(j) > s(j) - k THEN j := j+1 ELSE
        IF s(i) < s(j) THEN  i := j, j := j+1 ELSE
        A := A union {i}, i := j, j := j+1, state := up
\end{verbatim}
}
In other words, when it is time for an upward step, if the next value
goes up but not by $k$ ignore it, if it goes up by $k$ or more,
accept it as the new provisional value, and if it goes down,
replace the old provisional down step by the new value.
The pointer $i$ points to the provisional value at any time, and
when a new provisional value is accepted (rather than replacing
and old one), the old one becomes permanent.

\begin{lem} \label{lem:alg}
Let $s \in \cS_n$.  Then the subsequence $(s (i) : i \in A)$ defined
by the foregoing algorithm is a $k$-alternating subsequence of maximal
length.
\end{lem}

\noindent{\sc Proof:} Regarding $s$ as a word of length $n$, let
$s'$ denote the word of length $n-1$ obtained by removing 
the initial element of $s$ and let $s''$ denote the word of length $n-1$
obtained by removing the second element of $s$.  Let $L_{n,k}^*$
denote the length of the longest $k$-alternating sequence beginning
with a down step instead of an up step.  We claim that
\begin{eqnarray*}
s(1) < s(2) < s(1) + k & \Rightarrow & L_{n,k} (s) = L_{n,k} (s'') \\
s(1) > s(2) & \Rightarrow & L_{n,k} (s) = L_{n,k} (s') \\
s(1) + k \leq s(2) & \Rightarrow & L_{n,k} (s) = 1 +  L_{n,k}^* (s')
\end{eqnarray*}
The first holds because we can't use both $s(1)$ and $s(2)$ and starting
with $s_1$ dominates starting with~$s(2)$.  The second holds because
again we can't use both and this time starting from $s_2$ dominates
starting from~$s(1)$.  The last is true for the following reason.
The LHS cannot be more than the RHS because any $k$-alternating
subsequence restricts to a \emph{reverse} $k$-alternating sequence of
$s'$ upon removal of its first element (here the inequalities in the 
definition of alternating sequence are reversed, not the word itself).  
On the other hand, if $w$ is a reverse $k$-alternating subsequence of~$s'$, 
then there are two cases.  If the first element $w(1)$ is at least~$s(2)$, 
we can prepend $s(1)$ and obtain
a $k$-alternating subsequence of $s$ longer by one.  Similarly, if
the first element in $w$ is less than~$s(2)$, we can replace $w(1)$ 
by $s(2)$ and then prepend~$s(1)$.  This proves the claim.  The lemma
now follows by induction.
$\Cox$

Replacing $k$-alternation by $x$-alternation, an identical argument
shows that greedy provisional acceptance will also identify an
$x$-alternating subsequence of $\yy$ having maximal length.
Next, we adjust the bookkeeping slightly as follows.  The way
the algorithm is written, the first element $y_1$ begins in a state
of provisional acceptance.  When $y_1 > 1-x$, it is doomed eventually
to be replaced, so instead of provisionally accepting it, we reject
each initial value until we see a value that is at most $1-x$.
This yields the following easy lemma.
\begin{lem} \label{lem:x}
Conditional on $y_1 , \ldots , y_j$, the probability
of rejecting $y_{j+1}$ is always precisely $x$.
\end{lem}

\noindent{\sc Proof:} If no value has yet been provisionally accepted,
then by rule we reject precisely those values above $1-x$.  On the
other hand, if any value has been provisionally accepted, it is easy
to check inductively that when the state is ``up'', the provisional
value $y$ is at most $1-x$, and the rejection interval for the new
value, $[y,y+x)$ is entirely within $[0,1]$ and has length $x$.
Similarly, when the state is ``down'', the provisionally
accepted value is at least $x$ and the rejection interval
$(y-x,y]$ again has length $x$.
$\Cox$

Let $A \subseteq \{1,\ldots,n\}$ be the subset of indices $i$ for which
$y_i$ was at least provisionally accepted.  The previous
lemma shows that $A$ has the distribution of a set selected
by independent coin flips with success probability $1-x$.

\begin{lem} \label{lem:thinning}
Let $j_1 < j_2 < \cdots < j_r$ enumerate the set $A$.  Let
$z_i := y_{j_i}$ when $y_{j_i}$ was provisionally accepted
initially or after a down step and let $z_i := y_{j_i} - x$
when $y_i$ was provisionally accepted after an up step.
Then $\zz$ is a collection of independent variables uniform
on $[0,1-x]$ and is independent of $A$.
\end{lem}

\noindent{\sc Proof:}
Condition on the $y_1, \ldots , y_j$.  We know that
$\P (j+1 \in A) = 1-x$.  We therefore need to show that
conditional on $j+1 \in A$, and on $y_1, \ldots y_j$,
the value $z_{i+1}$ is uniform on $[0,1-x]$ where $i$
is the cardinality of $A \cap \{1,\ldots,j\}$.  When $i=0$ we are
in the initial phase and the result is obvious.  If not,
suppose first that the state is ``up''.  Then $z_i \leq 1-x$
and the values of $y_{j+1}$ for which provisional acceptance
will occur are the union of two intervals $[0,z_i] \cup [z_i+x,1]$.
If $y_{j+1}$ lies in the upper of these two intervals, it will
be provisionally accepted after an up step while if it is in
the lower interval it will be provisionally accepted replacing
a previous down step value.  Thus the two intervals together will
map to the single interval $[0,1-x]$.  Similarly, supposing
instead that the state is ``down'', provisional acceptance will
occur in $[0,z_i - x] \cup [z_i,1]$; $z_{i+1}$ will be $y_{j+1} - x$
in the upper interval and $y_{j+1}$ in the lower interval, and
again we see that $z_{i+1}$ is uniform on $[0,1-x]$.
$\Cox$

\noindent{\sc Proof of Lemma}~\ref{lem:select}:
Let $\zz$ be as in Lemma~\ref{lem:thinning}.
By Lemma~\ref{lem:alg}, the quantity $L_{n,x} (\yy)$ is equal
to $L_{|\zz|,0} (\zz)$.  By Lemma~\ref{lem:thinning}, the joint
distribution of $\disp{ \left ( |\zz| , \frac{\zz}{1-x} \right )}$
is the product measure ${\rm Bin}(n , 1-x) \times \mu$.  The
permutation associated with $\zz$ is the same as that associated
with the dilation $\zz / (1-x)$, whence the conditional distribution
of $\Psi (\zz)$ given $|\zz|$ is uniform on $\cS_{|\zz|}$, which is
enough fo prove the lemma.
$\Cox$

\subsection*{Final remarks}

The maximum of $(1-x)(2+5x)$ on $[0,1]$ occurs at $x = 3/10$.
Consequently the variance of the length of the longest $x$-alternating
sequence is maximized not at ordinary alternating sequences ($x=0$)
but at $0.3$-alternating sequences.

The asymptotics in Lemma~\ref{lem:compare} can be sharpened.  Instead
of tightness of the maximum of a Brownian bridge, use tightness of
the renomralized bridge statistic
$$\max \{ X(t) / \sqrt{t(1-t)} |\log (t(1-t))| \, : \, 0 \leq t \leq 1 \}.$$
This allows us to replace $x_2$ by $k/n + \min \{ n^{-1/3} ,
C (n-k)^{1/2 + \ee} \}$ in Lemma~\ref{lem:compare}.  The estimate
in Theorem~\ref{th:main} then becomes a sharp asymptotic $\E L_{n,k}
\sim (2/3) (n-k)$, uniform down to $n-k > (n-k)^{\delta}$,
where $\delta$ can be made arbitrarily small.


\end{document}